# EFFICIENT 1-BIT TENSOR APPROXIMATIONS

ALEX W. N. RIASANOVSKY AND SARAH EL KAZDADI

ABSTRACT. We present a spatially efficient decomposition of matrices and arbitrary-order tensors as linear combinations of tensor products of $\{-1, 1\}$-valued vectors. For any matrix $A \in \mathbb{R}^{m \times n}$,

$$A - R_w = S_w C_w T_w^\top = \sum_{j=1}^w c_j \cdot \boldsymbol{s}_j \boldsymbol{t}_j^\top$$

is a *w-width signed cut decomposition of A*. Here $C_w = \text{diag}(\boldsymbol{c}_w)$ for some $\boldsymbol{c}_w \in \mathbb{R}^w$, and $S_w, T_w$, and the vectors $\boldsymbol{s}_j, \boldsymbol{t}_j$ are $\{-1, 1\}$-valued.

To store $(S_w, T_w, C_w)$, we may pack $w \cdot (m + n)$ bits, and require only $w$ floating (or fixed-) point numbers. As a function of $w$, $\|R_w\|_F$ exhibits exponential decay when applied to `f32` matrices with i.i.d. $\mathcal{N}(0, 1)$ entries. Choosing $w$ so that $(S_w, T_w, C_w)$ has the same memory footprint as a `f16` or `bf16` matrix, the relative error is comparable.

Our algorithm yields efficient signed cut decompositions in 20 lines of pseudocode. It reflects a simple modification from a celebrated 1999 paper [1] of Frieze and Kannan.

As a first application, we approximate the weight matrices in the open `Mistral-7B-v0.1` Large Language Model to a 50% spatial compression. Remarkably, all 226 remainder matrices have a relative error $< 6\%$ and the expanded model closely matches `Mistral-7B-v0.1` on the `HuggingFace` leaderboard [2]. Benchmark performance degrades slowly as we reduce the spatial compression from 50% to 25%.

We optimize our open source `Rust` implementation [3] with `SIMD` instructions on `avx2` and `avx512` architectures. We also extend our algorithm from matrices to tensors of arbitrary order and use it to compress a picture of the first author's cat Angus.

## 1. INTRODUCTION

In a recent survey [4], Murray et al offer Randomized Numerical Linear Algebra (RNLA) as a treatment for large-scale problems in high performance computing (HPC) and machine learning (ML).

> *A dire situation.* While communities that rely on NLA now vary widely, they share one essential property: a ravenous appetite for solving larger and larger problems.
>
> — Murray et al [4]

Throughout the survey, the authors emphasize the emergence of structure at scale, which randomized algorithms are uniquely adept at exploiting.

The story in combinatorial limit theory is similar. Here, the *cut norm* $\|\cdot\|_\square$, first coined in 1999 in [1] by Frieze and Kannan, plays a vital role in exposing emergent structure in large networks. Using the cut norm, they approximate matrices with linear combinations of rank-1 matrices with entries in $\{0, 1\}$:

$$\|A\|_\square := \max_{S,T} \left| \sum_{(i,j) \in S \times T} a_{ij} \right| \rightsquigarrow A \approx \sum_j d_j \cdot \mathbf{1}_{S_j \times T_j}.$$

Following [1], researchers at the Theory Group of Microsoft Research applied these so-called *Frieze-Kannan decompositions* to connect several central, but superficially different, topics in combinatorics (c.f. [5])

Although Frieze-Kannan decompositions are spatially efficient, it is challenging to efficiently minimize their errors (c.f. Section 2.1). As a remedy, we base our methods around the signed cut norm $\|\cdot\|_\blacksquare$ instead of $\|\cdot\|_\square$, and instead approximate with outer products of $\{-1, 1\}$-valued vectors:

$$\|A\|_\blacksquare := \max_{\boldsymbol{s}, \boldsymbol{t}} \langle \boldsymbol{s} A, \boldsymbol{t} \rangle \rightsquigarrow A \approx \sum_j c_j \boldsymbol{s}_j \boldsymbol{t}_j^\top.$$

Summarizing Section 4, we argue that these *signed cut decompositions* produce state-of-the-art approximations. As noted in Section 3, the main algorithms we use to find them fit in 20 lines of pseudocode.

1.1. **Comparison to `bf16` quantization.** Our signed cut decompositions trade accuracy for space at a better exchange rate than `bf16` quantization offers for larger matrices, as shown in Figure 1. Moreover, we may make this tradeoff dynamically at runtime to suit the needs of our application.

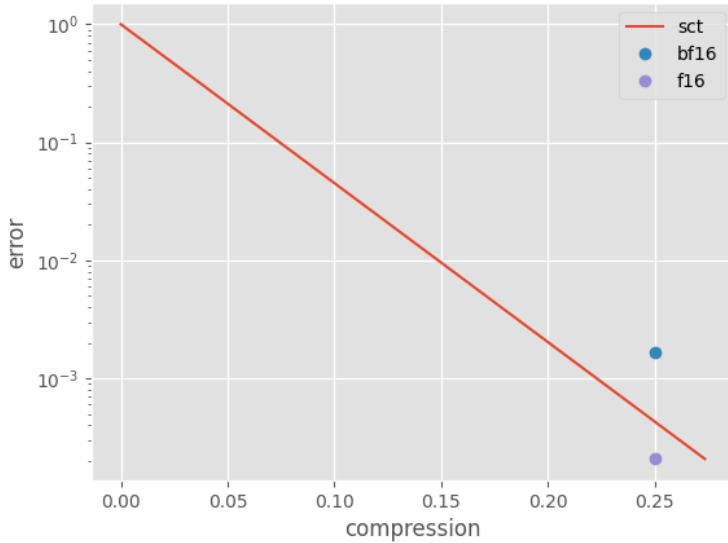

Figure 1. Relative error (log-scaled) and relative size of a signed cut decompositions of a $4096 \times 4096$ `f64` matrix with i.i.d. standard normal entries. For comparison, we note the tradeoff made with `f16` and `bf16`.

1.2. **Approximating `Mistral-7B-v0.1`.** When the original matrix is more structured, signed cut decompositions converge faster. In Figure 2, we show the effect of signed cut decompositions on the open `Mistral-7B-v0.1` model's benchmark performance.

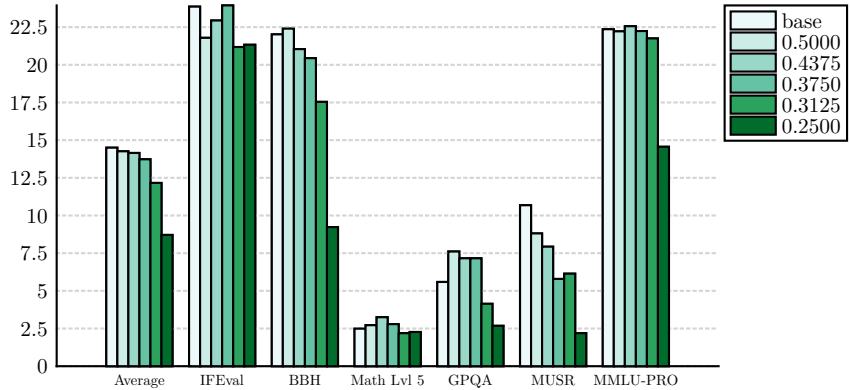

FIGURE 2. Benchmark results formed by approximating the weight matrices of `Mistral-7B-v0.1` with signed cut decompositions and expanding them at different widths. Here, 0.500 reflects a 2-fold spatial compression, and 0.2500 reflects a 4-fold spatial compression.

1.3. **Higher order tensors.** All of our techniques extend from matrices to tensors of arbitrary order, as show in Section 5. In Figure 3, we compress a picture the first author's cat by treating the implicit array of `RGB` values as $4000 \times 3000 \times 3$ tensor of real numbers.

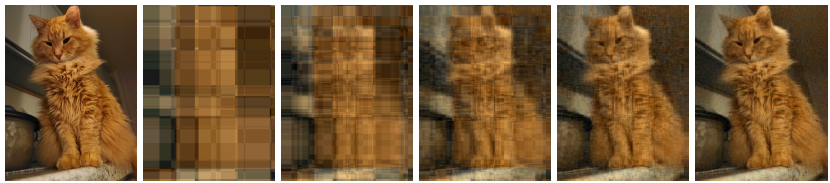

FIGURE 3. A picture of the first author's cat Angus and 5 approximations based on signed cut decompositions. The original `jpg` image occupies 792.2 `KiB` of memory. By comparison, an uncompressed `RGB` array of bytes would occupy 34.33 `MiB`. The selected approximations double in size from left (6.93 `KiB`) to right (110.88 `KiB`).

2. Cut norms and decompositions

2.1. **The unsigned cut norm and Frieze-Kannan decompositions.** We adapt techniques from combinatorial limit theory, a subfield of combinatorics cultivated by the Theory Group of Microsoft Research in the early 2000s. The germ of this theory is the *cut norm*, coined originally by Frieze and Kannan [1], which connects cut-set problems, matrix approximations, and the celebrated Regularity Lemma of Szemerédi. For any matrix $A = (a_{ij}) \in \mathbb{R}^{m \times n}$, the *cut norm* of $A$ is

$$\|A\|_\square := \max_{S,T} \left| \sum_{(i,j) \in S \times T} a_{ij} \right| = \max_{S,T} |\langle \mathbf{1}_{S \times T}, A \rangle|,$$

where $S$ ($T$) is any set of rows (columns) of $A$. In other words, $\|\cdot\|_\square$ measures the largest absolute sum over all submatrices of $A$. Although $\|A\|_\square$ is `MAX SNP hard`, it can theoretically be approximated efficiently in polynomial time (c.f. [6], [7]). As a motivating question, suppose we wish to optimize Frieze-Kannan decompositions on error, space, and time. As a function of the width $w$, we define the program

$$\text{FK}_w(A) : \begin{cases} \min & \|A - \sum_{j=1}^w d_j \cdot \mathbf{1}_{S_j \times T_j}\|_F \\ \text{s.t.} & d_1, ..., d_w \in \mathbb{R} \\ & S_1, ..., S_w \subseteq \{1, ..., m\} \\ & T_1, ..., T_w \subseteq \{1, ..., n\} \end{cases}.$$

In [1], the authors approach $\text{FK}_w(A)$ with randomized algorithms which estimate the cut norm of a residual matrix and use the result to perform a rank-1 update on the remainder matrix. Their algorithms are complicated by two unhelpful facts:

1. $\text{FK}_1(A)$ is not in general solved by calculating $\|A\|_\square$. The optimal approximation $A \approx d \cdot \mathbf{1}_{S \times T}$ has $d = \langle \mathbf{1}_{S \times T}, A \rangle \cdot (|S| \cdot |T|)^{-1}$ for some other $S, T$.
2. Pairs of matrices of the form $\mathbf{1}_{S \times T}$ are typically far from orthogonal. If we reuse our solution to $\text{FK}_w(A)$ to solve $\text{FK}_{w+1}(A)$, then these correlations force us in general to apply large updates to the scalars $d_1, ..., d_w$.

2.2. **Signed cut decompositions and the signed cut norm.** By instead basing our algorithms around the *signed cut norm* $\|\cdot\|_\blacksquare$, we mitigate the issues from Section 2.1 in practice. For any matrix $A \in \mathbb{R}^{m \times n}$, let

$$\|A\|_\blacksquare := \max_{s,t}\langle s, At\rangle = \max_{x,y}\langle x, Ay\rangle = \|A\|_{\infty \to 1}$$

where $s, t$ are $\{-1, 1\}$-valued and $x, y$ are $[-1, 1]$-valued. As show in [6], calculating $\|\cdot\|_\blacksquare$ is also `MAX SNP hard`. We approximate it with a randomized greedy algorithm. We define a *signed cut decomposition of width $w$* as any expression

$$A - R_w = S_w C_w T_w^\top = \sum_{j=1}^w c_j \cdot \mathbf{s}_j \mathbf{t}_j^\top$$

where $R_w \in \mathbb{R}^{m \times n}$, $S_w \in \{-1, 1\}^{m \times w}$, $C_w = \text{diag}(\mathbf{c}_w)$ for some $\mathbf{c}_w \in \mathbb{R}^w$, and $T_w \in \{-1, 1\}^{n \times w}$. As noted in Table 1, signed cut decompositions maintain the space efficiency of Frieze-Kannan decompositions while more closely emulating singular value decompositions (see Section 4).

| Decomposition | Rank-1 terms | Vector entries | Size (bits) |
|---|---|---|---|
| Singular value | $\sigma_j \cdot \mathbf{u}_j \mathbf{v}_j^\top$ | $\mathbb{R}$ | $(m + n + 1) \cdot f$ |
| Frieze-Kannan | $d_j \cdot \mathbf{1}_{S_j \times T_j}$ | $\{0, 1\}$ | $m + n + f$ |
| Signed cut | $c_j \cdot \mathbf{s}_j \mathbf{t}_j^\top$ | $\{-1, 1\}$ | $m + n + f$ |

TABLE 1. A comparison of three matrix decompositions equivalent to finite sums of rank-1 matrices. Here, a floating point number occupies $f$ bits in memory.

3. METHODS

For convenience, we let $\sigma(n) := \{-1, 1\}^n$. Our analog to $\text{FK}_w(A)$ is the program

$$\text{SC}_w(A) : \begin{cases} \min & \|A - \sum_{j=1}^{w} c_j \cdot \boldsymbol{s}_j \boldsymbol{t}_j^\top\|_F \\ \text{s.t.} & c_1, ..., c_w \in \mathbb{R} \\ & \boldsymbol{s}_1, ..., \boldsymbol{s}_w \in \sigma(m) \\ & \boldsymbol{t}_1, ..., \boldsymbol{t}_w \in \sigma(n) \end{cases}.$$

Our two main algorithms are straightforward and greedy. We include a few optimizations within each subsection and defer others (e.g., `SIMD`, bitset storage, and alignment) to Section 3.4. Algorithm 4 estimates $\|\cdot\|_\blacksquare$ iteratively. Algorithm 5 bootstraps these estimates to extend approximate solutions from $\text{SC}_k(A)$ to $\text{SC}_{k+1}(A)$. Finally, Algorithm 6 adjusts the coefficients to improve the regression. In Section 5, we generalize to tensors.

3.1. **Naive signed cut-sets.** Note that $\|A\|_\blacksquare$ is the solution to the program

$$\text{SgnCut}(A) : \begin{cases} \max & \langle \boldsymbol{s}, A\boldsymbol{t} \rangle \\ \text{s.t.} & \boldsymbol{s} \in \sigma(m) \\ & \boldsymbol{t} \in \sigma(n) \end{cases}.$$

In [1], the authors explore the analogous program for the cut norm. By replacing $\boldsymbol{s}$ with $\text{sgn}(A\boldsymbol{t})$, or $\boldsymbol{t}$ with $\text{sgn}(A^\top \boldsymbol{s})$, we see that $\text{Cut}(A)$ is equivalent to the programs

$$\text{LeftSgnCut}(A) : \begin{cases} \max & \|A^\top \boldsymbol{s}\|_1 \\ \text{s.t.} & \boldsymbol{s} \in \sigma(m) \end{cases} \qquad \text{RightSgnCut}(A) : \begin{cases} \max & \|A\boldsymbol{t}\|_1 \\ \text{s.t.} & \boldsymbol{t} \in \sigma(n) \end{cases}.$$

Algorithm 4 exploits this to approach a locally maximal cut by flipping signs.

---

<u>GREEDY SIGNED CUT</u>$(A \in \mathbb{R}^{m \times n})$:

1   **sample** $\boldsymbol{s}_0 \in \sigma(m)$ uniformly
2   **sample** $\boldsymbol{t}_0 \in \sigma(n)$ uniformly
3   **let** $c_0 := -\infty$
4   **for** $j \in \{1, 2, ...\}$**:**
5       **let** $\boldsymbol{s}_j := \text{sgn}(A \cdot \boldsymbol{t}_{j-1})$
6       **let** $\boldsymbol{t}_j := \text{sgn}(A^\top \cdot \boldsymbol{s}_j)$
7       **let** $c_j := \langle \boldsymbol{s}_j, A \cdot \boldsymbol{t}_j \rangle$
8       **if** $c_j \leq c_{j-1}$**:**
9           **yeet** $(c_j, \boldsymbol{s}_j, \boldsymbol{t}_j)$

ALGORITHM 4. A randomized greedy algorithm for $\text{Cut}(A)$.

By caching the intermediate vectors $A \cdot \boldsymbol{t}_j, A^\top \boldsymbol{s}_j$, we may sparsify the matrix-vector products. To see this, note that at iteration $j \in \{2, 3, ...\}$,

$$\boldsymbol{s}_j = \text{sgn}(A \cdot \boldsymbol{t}_{j-1}) = \text{sgn}(A \cdot \boldsymbol{t}_{j-2} + A \cdot (\boldsymbol{t}_{j-1} - \boldsymbol{t}_{j-2})), \text{ and}$$
$$\boldsymbol{t}_j = \text{sgn}(A^\top \cdot \boldsymbol{s}_{j-1}) = \text{sgn}(A^\top \cdot \boldsymbol{s}_{j-2} + A^\top \cdot (\boldsymbol{s}_{j-1} - \boldsymbol{s}_{j-2})).$$

In practice, as $j$ increases, the $\{-2, 0, 2\}$-valued vectors $\boldsymbol{s}_j - \boldsymbol{s}_{j-1}$ and $\boldsymbol{t}_j - \boldsymbol{t}_{j-1}$ become increasing sparse. We store the variables as follows:

- **sign bitsets**: the $\{-1, 1\}$-valued vectors $\boldsymbol{s}_j$ and $\boldsymbol{t}_j$ are stored with unsigned integers. Here, each 1 bit indicating a $-1$. At iteration $j$, we also store $\boldsymbol{s}_{j-1}$ and $\boldsymbol{t}_{j-1}$.
- **contiguous, aligned vectors**: the real-valued vectors $A \cdot \boldsymbol{t}_{j-1}$ and $A^\top \cdot \boldsymbol{s}_j$ are stored with aligned slices of `f32` values. The matrices $A$ and $A^\top$ are aligned, padded, and in column-major.

3.2. **Simple signed cut decompositions.** Algorithm 5 builds signed cut decompositions by repeatedly invoking Algorithm 4.

---

$\underline{\text{GREEDY DECOMPOSITION}}(A \in \mathbb{R}^{m \times n}, w \in \mathbb{N})$:

1  **let** $R_0 := A$
2  **let** $S_0 := 0 \in \{-1, 1\}^{m \times 0}$
3  **let** $\boldsymbol{c}_0 := 0 \in \mathbb{R}^0$
4  **let** $T_0 := 0 \in \{-1, 1\}^{n \times 0}$
5  **for** $k \in \{0, ..., w-1\}$:
6      **let** $(c_{k+1}, \boldsymbol{s}_{k+1}, \boldsymbol{t}_{k+1}) = \texttt{greedy\_signed\_cut}(R_k)$
7      **let** $\boldsymbol{c}_{k+1} := [\boldsymbol{c}_k \ c_{k+1}]$
8      **let** $S_{k+1} := [S_k \ \boldsymbol{s}_{k+1}]$
9      **let** $T_{k+1} := [T_k \ \boldsymbol{t}_{k+1}]$
10     **let** $R_{k+1} := R_k - \frac{c_{k+1}}{mn} \cdot \boldsymbol{s}_{k+1} \boldsymbol{t}_{k+1}^\top$
11 **yeet** $(S_w, \mathrm{diag}(\boldsymbol{c}_w), T_w)$

---

ALGORITHM 5. A "greedy" least-squares approach to signed cut decompositions. The remainder $R_{k+1}$ is the orthogonal complement of $R_k$ with respect to the 1-dimensional subspace of $\mathbb{R}^{m \times n}$ spanned by $\boldsymbol{s}_{k+1} \boldsymbol{t}_{k+1}^\top$.

By inspection,

$$\|R_{k+1}\|_F^2 = \|R_k\|_F^2 - \frac{c_{k+1}^2}{mn}.$$

for all $0 \leq k < w$. If $c_k$ is proportional to $\|R_{k-1}\|_F$, then $\|R_k\|_F$ decreases exponentially. Empirically, results in Section 4 affirm this hypothesis.

To amortizes the cost of the rank-1 updates to $R_k$, we delay the `matmul` accumulation in Algorithm 5 by storing each $R_k$ implicitly as $(R', S, C, T)$ where $R_k = R' + SCT^\top$. When $(S, C, T)$ reaches a carefully chosen width (32), we flush $SCT^\top$ into $R'$ to form $R_k$. We store all real-valued vectors and matrices (except the diagonal matrices) as in Section 3.1. For the $\{-1, 1\}$-valued matrices, we also use bitsets with the 1 bit indicating a $-1$ sign.

3.3. **Least square corrections.**

There is a simple improvement that can be made to Algorithm 5 based on the method of least squares. In the $k$-th iteration of Algorithm 6, we treat $S_k$ and $T_k$ as fixed and allow the diagonal matrix $C_k$ to vary so as to minimize $\|A - S_k C_k T_k^\top\|_F$.

---

LEAST SQUARES DECOMPOSITION($A \in \mathbb{R}^{m \times n}$, $w \in \mathbb{N}$):

1  **let** $S_0 := 0 \in \{-1, 1\}^{m \times 0}$
2  **let** $\boldsymbol{c}_0 := 0 \in \mathbb{R}^0$
3  **let** $T_0 := 0 \in \{-1, 1\}^{n \times 0}$
4  **for** $k \in \{0, ..., w-1\}$**:**
5      **let** $R_k := S_k \, \mathrm{diag}(\boldsymbol{c}_k) T_k^\top$
6      **let** $(\_, \boldsymbol{s}_{k+1}, \boldsymbol{t}_{k+1}) = \mathtt{greedy\_signed\_cut}(R_k)$
7      **let** $S_{k+1} := [S_k \; \boldsymbol{s}_{k+1}]$
8      **let** $T_{k+1} := [T_k \; \boldsymbol{t}_{k+1}]$
9      **let** $\boldsymbol{c}_{k+1} := \mathrm{argmin} \, \|A - S_{k+1} \, \mathrm{diag}(\boldsymbol{c}_{k+1}) T_{k+1}^\top\|_F$
10 **yeet** $(S_w, \mathrm{diag}(\boldsymbol{c}_w), T_w)$

---

ALGORITHM 6. An improvement to Algorithm 5. Here, $R_{k+1}$ is the orthogonal complement of $A$ with respect to the span of $\boldsymbol{s}_1 \boldsymbol{t}_1^\top, ..., \boldsymbol{s}_{k+1} \boldsymbol{t}_{k+1}^\top$.

To calculate $\boldsymbol{c}_{k+1}$ in Algorithm 6, we formulate the associated least-squares problem and solve its *normal equation*

$$X_k^\top X_{k+1} \cdot \boldsymbol{c}_{k+1} = X_{k+1}^\top \cdot A.$$

Here, $X_{k+1}$ is the linear operator from $\mathbb{R}^k$ to $\mathbb{R}^{m \times n}$ defined by sending $j$-th standard basis vector to $\boldsymbol{s}_j \boldsymbol{t}_j^\top$.

In our experiments, the improvement from swapping Algorithm 5 for Algorithm 6 is not significant for many matrices. We offer two possible explanations:

1. If only $c_{k+1}$ varies, then Algorithm 5 and Algorithm 6 are equivalent.
2. The more orthogonal the matrices $\boldsymbol{s}_k \boldsymbol{t}_k^\top$ are, the less Algorithm 6 helps.

A lengthier analysis goes beyond the scope of this paper.

3.4. **SIMD Optimizations.** We improve the implementation in [3] by leveraging `SIMD` intrinsics on both `avx2` and `avx512` architectures. Interestingly, neither Algorithm 4 nor Algorithm 5 makes a single call to `fused multiply add` instructions. To demonstrate this, suppose we wish to compute the inner product $\langle \boldsymbol{s}, \boldsymbol{x} \rangle$ where

$$\boldsymbol{s} = \begin{bmatrix} s_1 \\ s_2 \\ \vdots \\ s_n \end{bmatrix} \in \sigma(n) \quad \text{and} \quad \boldsymbol{x} = \begin{bmatrix} x_1 \\ x_2 \\ \vdots \\ x_n \end{bmatrix} \in \mathbb{R}^n \quad .$$

Assuming $n = 8\ell$ for some $\ell$, we apply the isometry $\mathbb{R}^n \cong \bigoplus_{j=1}^\ell \mathbb{R}^8$ and write

$$\boldsymbol{s} = \boldsymbol{s}_1 \oplus \cdots \oplus \boldsymbol{s}_\ell \qquad \text{where each } \boldsymbol{s}_j \in \sigma(8), \text{ and}$$
$$\boldsymbol{x} = \boldsymbol{x}_1 \oplus \cdots \oplus \boldsymbol{x}_\ell \qquad \text{where each } \boldsymbol{x}_j \in \mathbb{R}^8, \text{ so that}$$
$$\langle \boldsymbol{s}, \boldsymbol{x} \rangle = \langle \boldsymbol{s}_1, \boldsymbol{x}_1 \rangle + \cdots + \langle \boldsymbol{s}_\ell, \boldsymbol{x}_\ell \rangle.$$

Since our signed sets are stored as bitsets, we need only use the bits of each $\boldsymbol{s}_j$ to flip the sign of each $\boldsymbol{x}_j$ entry and then sum. On `avx512` instructions (and the upcoming `avx10` instructions), we may offload this task to `mask` instructions with 512-bit registers. On `avx2` and `aarch64`, we may instead broadcast 32 bits to smaller `f32` registers, then use shifts instructions to move bits of interest to flip sign bits.

## 4. Results

4.1. **Approximating a random matrix.** In this section, we evaluate the quality of signed cut decompositions by approximating a random matrix $A \in \mathbb{R}^{4096 \times 4096}$ with independent $\mathcal{N}(0,1)$ entries. As an array of `f64` numbers, this occupies $8 \cdot 4096 \cdot 4096$ bytes, or $64 \cdot 4096 \cdot 4096$ bits, in memory. In Figure 7, we consider the intermediate signed cut decompositions $A = R_k + \sum_{j=1}^{k} c_j \boldsymbol{s}_j \boldsymbol{t}_j^\top$. We store each triple $(c_j, \boldsymbol{s}_j, \boldsymbol{t}_j)$ so the scalar has `f64` precision and the signed vectors $\boldsymbol{s}_j, \boldsymbol{t}_j$ as a bitset, as in Section 3. This way, each triple can be stored with $64 + m + n$ bits in memory. For each $k$, we plot $(p_k, r_k)$ where

$$p_k := \frac{k \cdot (64 + 4096 + 4096)}{64 \cdot 4096 \cdot 4096} \quad \text{and} \quad r_k := \frac{\|R_k\|_F}{\|A\|_F}.$$

This way, $p_k$ is the *compression rate* of the $k$-width approximation of $A$. For comparison, we plot $(0.25, r)$ where $r$ is the relative error from downcasting $A$ to `bf16` or `f16` precision.

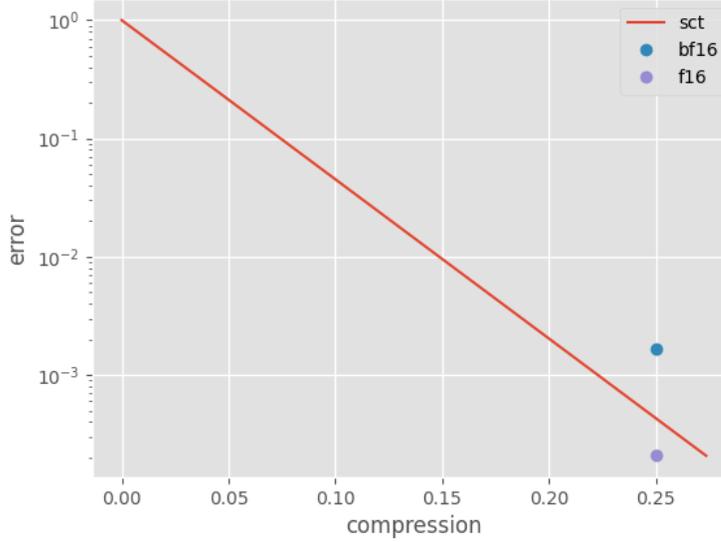

...EFFICIENT 1-BIT TENSOR APPROXIMATIONS 9FIGURE 7. Relative error (log-scale) plotted against spatial compression rate for a standard normal $4096 \times 4096$ matrix with `f64` precision. Here, `bf16` has error $\approx 0.00166$ and `f16` has error $\approx 0.000208$. Signed cut decompositions breaks even on error at $0.2064$ and $0.2734$, respectively.

## 4.2. Approximating `Mistral-7B-v0.1`.

In this section, we apply signed cut decompositions to approximate the open `Mistral-7B-v0.1` model.

4.2.1. *Frobenius norm error.* We repeat our experiment in Section 4.1 by forming signed cut decompositions of the 226 weight matrices from the open `Mistral-7B-v0.1` model. Since these matrices are stored with `bf16` precision, each such $A \in \mathbb{R}^{m \times n}$ occupies $16mn$ bits in the original `safetensors` file. In the approximation $A \approx \sum_{j=1}^{k} c_j \boldsymbol{s}_j \boldsymbol{t}_j^\top$, we store each scalar $c_j$ with `f32` precision and store each vector sign vector $\boldsymbol{s}_j, \boldsymbol{t}_j$ as a bitset, as described in Section 3. So the summands together occupy $32 + m + n$ bits in memory.

Targeting a 2-fold spatial compression, our approximation $A \approx \sum_{j=1}^{w} c_j \boldsymbol{s}_j \boldsymbol{t}_j^\top$ has width

$$w \approx \frac{1}{2} \cdot \frac{16mn}{32 + m + n}.$$

In Figure 8, we plot the relative error as a function of the compression rate.

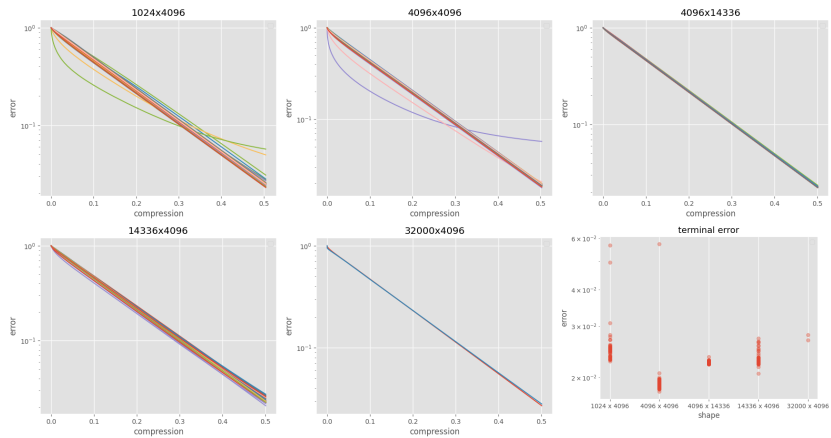

FIGURE 8. Relative error $\frac{\|R_k\|_F}{\|A\|_F}$ (log-scale) of signed cut decompositions of every matrix in `Mistral-7b-v0.1` after running Algorithm 5. Here, we group by matrix shape. The final width corresponds to a 50% spatial compression from the original `bf16` matrix.

4.2.2. *Large language model benchmarks.* Using $\frac{j}{16}$ for $j \in \{4, 5, 6, 7, 8\}$ as our target compression rate, we select appropriate widths based on the matrix dimensions $m, n$ as in Table 2.

| Shape | Count | Width by compression | | | | |
|---|---|---|---|---|---|---|
| | | 0.5000 | 0.4375 | 0.3750 | 0.3125 | 0.2500 |
| $1024 \times 4096$ | 64 | 6512 | 5698 | 4884 | 4070 | 3256 |

| Shape | Count | Width by compression | | | | |
|---|---|---|---|---|---|---|
| | | 0.5000 | 0.4375 | 0.3750 | 0.3125 | 0.2500 |
| $4096 \times 4096$ | 64 | 16320 | 14280 | 12240 | 10200 | 8160 |
| $4096 \times 14336$ | 32 | 25442 | 22261 | 19081 | 15901 | 12721 |
| $14336 \times 4096$ | 64 | 25442 | 22261 | 19081 | 15901 | 12721 |
| $32000 \times 4096$ | 2 | 29023 | 25395 | 21767 | 18139 | 14511 |

Table 2. Widths matching each shape to each desired compression rate.

We expand the corresponding truncated signed cut approximations and write new safetensors files, reusing all other values. In order to measure the quality of signed cut decompositions as a vehicle for model quantization, we truncate the weight matrix approximations $A \approx \sum_{j=1}^{w} c_j s_j t_j^\top$ at different weights corresponding to our target compression rates. With $m, n$ fixed, a $w$-width signed cut approximation of at bf16 $m \times n$ matrix has compression rate

$$p_{m,n}(w) = \frac{w \cdot (32 + m + n)}{16mn}.$$

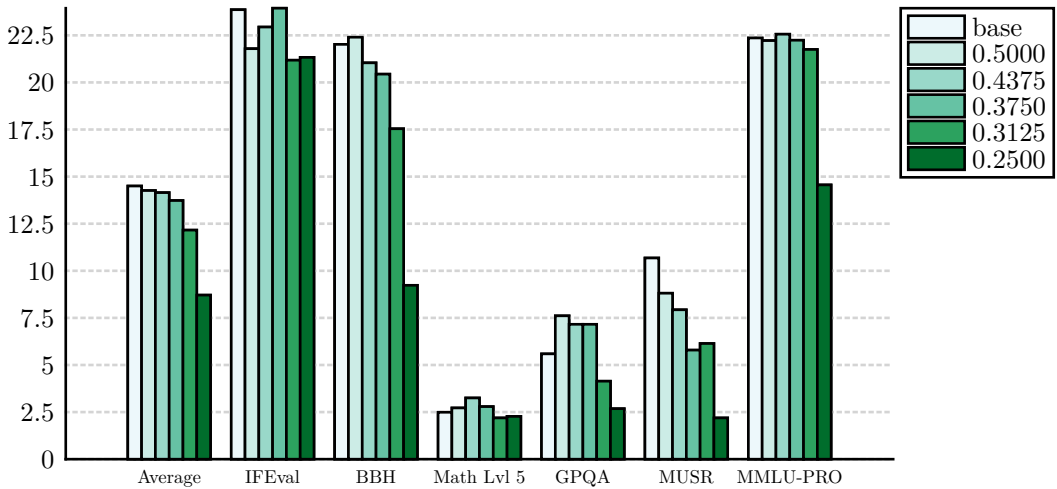

Figure 9. HuggingFace benchmark [2] performance of signed cut decompositions of the open Mistral-7B-v0.1 model at compression rates between 0.5 and 0.25.

## 5. Generalizations

As mentioned in Section 1, we motivate the tensor signed cut decomposition by approximating a picture of the first author's cat Angus, as in Figure 10. To begin, we convert this jpg image to a $4000 \times 3000$ array of RGB triples, with each value stored as an unsigned 8-bit integer. Diving along the "color" axis, we have 3 $4000 \times 3000$ matrices corresponding to the colors red, green, and blue. To exploit correlations among these 3 matrices, we find signed cut approximations of each using a common sequences $s_1, \ldots$ and $t_1, \ldots$ of sign vectors. Our goal is to approx-imate image tensors of the form $\boldsymbol{a} \in \mathbb{R}^{m \times n \times 3}$ where $\min(m,n) \gg 3$. We compare two variants of this decomposition:

1. **RGB signs**: approximate with tensors of the form $\sum_j c_j \cdot \boldsymbol{s}_j \otimes \boldsymbol{t}_j \otimes \boldsymbol{k}_j$ where $\boldsymbol{s}_j$, $\boldsymbol{t}_j$, and $\boldsymbol{k}_j$ are $\{-1,1\}$-valued. Here, we minimize the number of bits used to store vectors but sacrifice some flexibility.
2. **RGB scalars**: approximate with tensors of the form $\sum_j \boldsymbol{s}_j \otimes \boldsymbol{t}_j \otimes \boldsymbol{c}_j$ where $\boldsymbol{s}_j$, and $\boldsymbol{t}_j$, and $\boldsymbol{c}_j \in \mathbb{R}^3$. As shown in Figure 10, this approach converges noticeably faster.

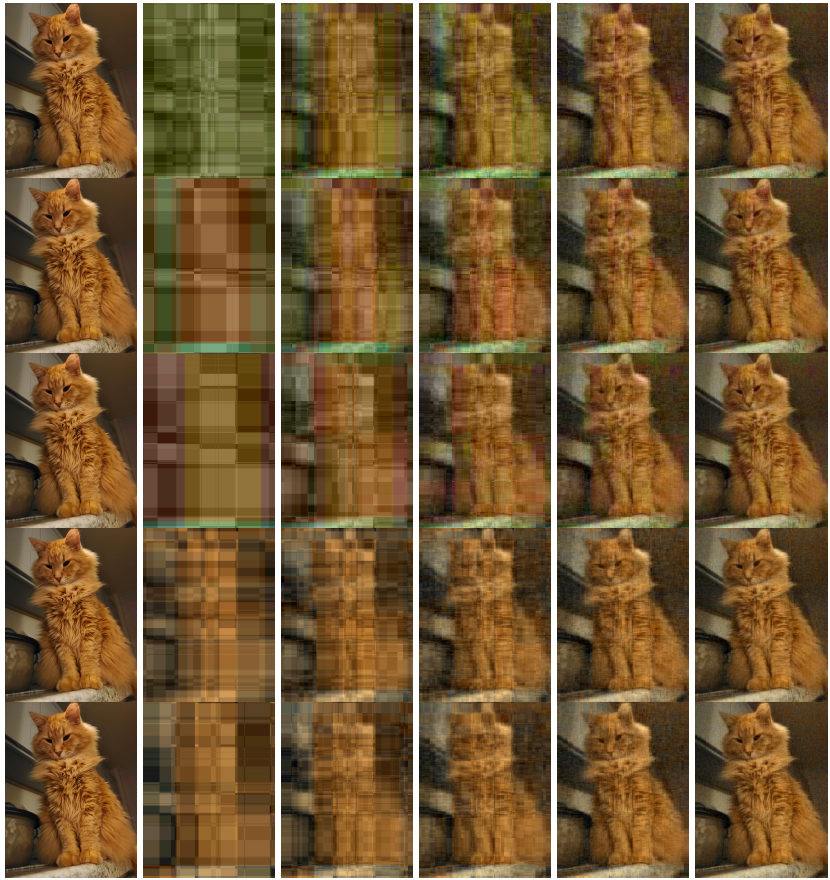

FIGURE 10. Approximations of the first author's cat Angus with different combinations of heuristics. The first 3 rows use RGB signs: (first with Algorithm 12, then with a least squares correction, and finally a costly heuristic to optimize sign vectors). On the remaining 2 rows, we use RGB scalars (first with least squares, then with costly sign vector optimizations). C.f. [3].

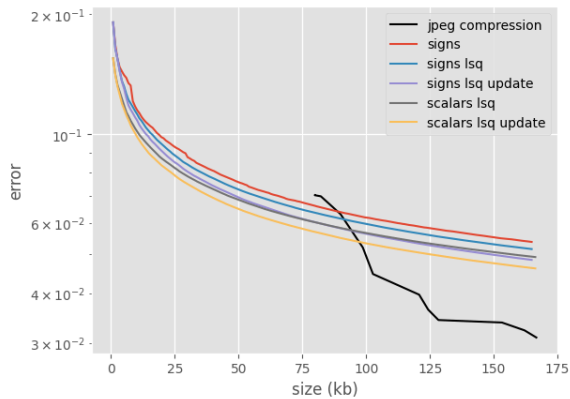

FIGURE 11. $L_2$ errors (log-scale) corresponding to the approximations in Figure 10. Here, we normalize error between tensors in $\{0, ..., 255\}^{m \times n \times 3}$ by $255 \cdot \sqrt{3mn}$. For comparison, we apply `jpg` compression to the original image with "quality" ranging from 1% to 13%.

5.1. **Signed cuts.** For full generality, we generalize our approach (namely, `RGB` signs from Section 5) to arbitrary higher-order tensors. For our purposes, an order-$k$ tensor is an array $\boldsymbol{a} \in \mathbb{R}^{\otimes \boldsymbol{n}}$ where $\boldsymbol{n} = (n_1, ..., n_k) \in \mathbb{N}^k$. We define the signed cut norm over $\mathbb{R}^{\otimes \boldsymbol{n}}$ as follows. First, let

$$\sigma(\boldsymbol{n}) := \{-1, 1\}^{\otimes \boldsymbol{n}} = \bigotimes_{i=1}^{k} \sigma(n_i).$$

Here, $\sigma(\boldsymbol{n})$ is a $2^{n_1 + \cdots + n_k}$-element frame spanning $\mathbb{R}^{\otimes \boldsymbol{n}}$. Then for all $\boldsymbol{a} \in \mathbb{R}^{\otimes \boldsymbol{n}}$,

$$\|\boldsymbol{a}\|_{\blacksquare} := \max_{\boldsymbol{s} \in \sigma(\boldsymbol{n})} \langle \boldsymbol{s}, \boldsymbol{a} \rangle.$$

Generalizing from Section 3.1, we note that for all $\boldsymbol{n} \in \mathbb{N}^k$ and all $1 \leq j \leq k$, the programs

$$\text{TensorSgnCut}(A) : \begin{cases} \max & \langle \boldsymbol{s}, \boldsymbol{a} \rangle \\ \text{s.t.} & \boldsymbol{s} \in \sigma(\boldsymbol{n}) \end{cases} \text{ and }$$

$$\text{AxialCut}_i(A) : \begin{cases} \max & \|\boldsymbol{a} \times_i \boldsymbol{s}\|_1 \\ \text{s.t.} & \boldsymbol{s} \in \sigma\left(\left(n_j\right)_{j \neq i}\right) \end{cases}$$

both solve for $\|\boldsymbol{a}\|_{\blacksquare}$. Here, $\boldsymbol{a} \times_i \boldsymbol{s} \in \mathbb{R}^{n_j}$ is the product of $\boldsymbol{a}$ and $\boldsymbol{s}$ along the $i$-th axis. We extend Algorithm 4 from Section 3.1 to Algorithm 12.

```
AXIAL GREEDY CUT(a ∈ ℝ^⊗n):
1   sample s_0 = ⊗_{i∈A} s_{0i} ∈ σ(n) uniformly
2   let c_0 := -∞
3   for j ∈ {1, 2, ...}:
4       for i ∈ {1, ..., k}:
5           let s_{j↓i} := ⊗_{i'<i} s_{ji}
6           let s_{j↑i} := ⊗_{i'>i} s_{j-1,i}
7           let s_{j↕i} := s_{j↓i} ⊗ s_{j↑i}
8           let s_{ji} := sgn(a ×_i s_{j↕i})
9       let s_j := ⊗_{i∈A} s_{ji}
10      let c_j := ⟨s_j, a⟩
11      if c_j ≤ c_{j-1}:
12          yeet (c_j, s_j)
```

ALGORITHM 12. Algorithm 4 for tensors.

**5.2. Signed cut decompositions.** Here, for any $a \in \mathbb{R}^{\otimes n}$, a *signed cut decomposition of width $w$* is any solution to

$$a - r_w = S_w \cdot c_w = \sum_{j=1}^{w} c_j \cdot s_j$$

where each $s_j \in \sigma(n)$, $S_w \in \sigma(n) \otimes \mathbb{R}^w$ can be treated as a linear operator from $\mathbb{R}^w$ to $\mathbb{R}^{\otimes n}$ which sends the $j$-th standard basis vector to $s_j$, $c_w = (c_j) \in \mathbb{R}^w$, and $r_w \in \mathbb{R}^{\otimes n}$. Continuing, we extend Algorithm 5 from Section 3.2 to Algorithm 13.

```
TENSOR GREEDY DECOMPOSITION(a ∈ ℝ^⊗n, w ∈ ℕ):
1   let r_0 := a
2   let S_0 := 0 ∈ σ(n) ⊗ ℝ^0
3   let c_0 := 0 ∈ ℝ^0
4   for k ∈ {0, ..., w-1}:
5       let (c_{k+1}, s_{k+1}) = axial_greedy_cut(r_k)
6       let c_{k+1} := [c_k  c_{k+1}]
7       let S_{k+1} := [S_k  s_{k+1}]
8       let r_{k+1} := r_k - (c_{k+1} / ∏_i n_i) · s_{k+1}
9   yeet (c_w, S_w)
```

ALGORITHM 13. Algorithm 5 for tensors.

5.3. **Least square corrections.** Finally, we extend Algorithm 6 from Section 3.3 to Algorithm 14.

---

$\underline{\text{Tensor least Squares}}(\boldsymbol{a} \in \mathbb{R}^{\otimes \boldsymbol{n}}, w \in \mathbb{N})$:

1   **let** $S_0 := 0 \in \sigma(\boldsymbol{n}) \otimes \mathbb{R}^0$

2   **let** $\boldsymbol{c}_0 := 0 \in \mathbb{R}^0$

3   **for** $k \in \{0, ..., w-1\}$**:**

4       **let** $\boldsymbol{r}_k := \boldsymbol{a}_k - S_k \cdot \boldsymbol{c}_k$

5       **let** $(\_, \boldsymbol{s}_{k+1}) = \texttt{axial\_greedy\_cut}(\boldsymbol{r}_k)$

6       **let** $\boldsymbol{c}_{k+1} := [\boldsymbol{c}_k \; \boldsymbol{c}_{k+1}]$

7       **let** $S_{k+1} := [S_k \; \boldsymbol{s}_{k+1}]$

8       **let** $\boldsymbol{c}_{k+1} := \arg\min \|\boldsymbol{a} - S_{k+1} \cdot \boldsymbol{c}_{k+1}\|_2$

9   **yeet** $(\boldsymbol{c}_w, S_w)$

ALGORITHM 14. Algorithm 6 for tensors.

6. Discussion and Future Work

We discuss some potential applications of signed cut decompositions for future work.

6.1. **High Performance Computing.** In Section 4.1, we showed that signed cut decompositions can approximate large random matrices roughly as well as `bf16` and `f16` quantization while matching the memory footprint. In Section 4.2, we also showed that they approximate highly structured low-precision matrices even after driving down the memory footprint by a factor of 2 or higher. How do these gains in space complexity yield better runtimes? More precisely, we ask:

**Question A**: How do signed cut decompositions perform against other methods of quantization for `matvec`, `matvec`, and other common tasks?

At this time, we lack the resources to answer this question as thoroughly as we would like. Based on our benchmarks in [3], we estimate that signed cut decompositions cost space and time within the same order of magnitude as other quantization approaches which achieve the same accuracy. Since signed cut decompositions can dynamically sacrifice accuracy by truncating approximations at lower widths, and the correct widths likely vary heavily across different domains.

In our experiments, the width of typical signed cut decompositions frequently should be much larger than the dimensions of the original matrix. Optimal performance depends heavily on cache efficiency, which may require specialized layouts based on the number of rows, columns, and desired width.

Additionally, as noted in Section 3.4, our `matmul` implementation flips bits and sum instead of taking traditional inner products and is unable to leverage optimized

`FMA` instructions. Based on our preliminary experiments with fixed point numbers, we are especially curious how they may performance with signed cut decompositions.

6.2. **Images Processing.** As shown with Figure 10 and Figure 11, signed cut decompositions provide relatively accurate approximations of images when space is heavily constrained. Moreover, these decompositions admit straightforward data-parallel layouts, unlike traditional schemes for image compression and approximation.

**Question B**: Do there exist image processing workloads which can leverage signed cut decompositions? Furthermore, can we improve spatial efficiency while maintaining data-parallelism by restricting $s_j$ and $t_j$ to highly structured vectors?

For example, it may be efficient to assume that signs flip infrequently, or periodically, and optimize storage schema accordingly.

6.3. **Machine Learning.** Recall from Section 4.2 that our experiments replace all 226 weight matrices from the open `Mistral-7B-v0.1` model with the result of expanding signed cut decompositions $\sum_{j=1}^{w} c_j s_j t_j^\top$ where the width $w$ is chosen so as to emulate a predetermined compression rate. Based on Table 2, we see that our 2-fold compression consists of

$$64 \cdot 6{,}512 + 64 \cdot 16{,}320 + 64 \cdot 25{,}442 + 32 \cdot 25{,}442 + 2 \cdot 29{,}023 = 3{,}961{,}726$$

scalars $c_j$ (and the same number of row and column vectors $s_j$ and $t_j$). In other words:

**Question C**: If the $s_j$ and $t_j$ are fixed, can the `Mistral-7B-v0.1` model be efficiently retrained treating only the 3.96 million scalars as variables?

This question goes far beyond our current resources but may be of interest to other researchers.

For our final question, we remark that any depth-1 feed-forward universal approximation schema

$$f(\boldsymbol{x}) \approx f_{\mathrm{NN}}(\boldsymbol{x}) = A_2 \cdot \sigma(\boldsymbol{b}_1 + A_1 \cdot \boldsymbol{x})$$

can be transformed to a universal approximation schema

$$f(\boldsymbol{x}) \approx g_{\mathrm{NN}}(\boldsymbol{x}) = L_2 C_2 \cdot \sigma(\boldsymbol{b}_1 + L_1 C_1 R_1 \cdot \boldsymbol{x})$$

where $L_1, R_1, L_2$ have entries in $\{-1, 1\}$ and $C_1, C_2$ are diagonal matrices.

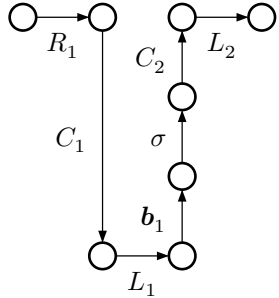

Figure 15. A universal approximation schema consisting of diagonal matrices $C_1, C_2$ and $\{-1, 1\}$-valued matrices $R_1, L_1, R_2$.

**Question D**: What is the optimal way to train using the architecture in Figure 15? What modifications to this architecture perform best?

With activation function $\sigma$ fixed, we are free to optimize over:

1. **shapes**: here, $g_{\text{NN}}$'s layers have dimensions $(m, k, \ell, n)$ where $k, \ell$ are free
2. **signs**: there are $nk + k\ell + \ell m$ variables comprising $L_1, R_1, R_2$
3. **reals**: $k + 2\ell$ real variables comprise $\boldsymbol{b}_1, C_1, C_2$

We leave this experiment as a topic for other researchers.

## 7. Acknowledgements

We thank Michael Sklar and the anonymous reviewers for providing helpful comments. Additionally, we thank Bernard Lidický for access to a high performance computing cluster which ran our first prototype.

Philadelphia, PA, USA
*Email address:* a.riasanovsky@gmail.com